\documentclass[a4paper,reqno]{amsart}
      \usepackage{etex}
      \usepackage[OT2, T1]{fontenc}
      \usepackage{url}
      \usepackage{amsmath}
      \usepackage{lmodern}
      \usepackage{array}
      \usepackage{graphicx}
      \usepackage{amsfonts}
      \usepackage{amssymb}
      \usepackage{extpfeil}
\usepackage{soul}
\usepackage{xargs}
      \usepackage{xcolor}
      \definecolor{imperialred}{RGB}{237, 41, 57}
      \definecolor{royalblue}{RGB}{64, 106, 212}
      \definecolor{link}{RGB}{11,0,128}
      \definecolor{olivegreen}{RGB}{128, 128, 0}
      \usepackage[colorlinks=true, citecolor=royalblue,linkcolor=olivegreen,urlcolor=link]{hyperref}
      \usepackage{amsthm}
      \usepackage{esint}
      \usepackage{epic}
      \usepackage{rotating}

      \usepackage{calligra}
      \usepackage{xkeyval}

      \usepackage{exscale, relsize}
      \usepackage{stackengine}
      \usepackage{scalerel}
      \usepackage{setspace}
      \usepackage{paralist}
      \usepackage{colonequals}		
      \usepackage{mathabx}		
      \usepackage[left]{lineno}
      \usepackage{amscd}
      \usepackage[all,cmtip]{xy}	
      \usepackage{sseq}			
      \usepackage{verbatim}
      \usepackage{parskip}
      \usepackage{microtype}
      \usepackage{ragged2e}
      \usepackage[in]{fullpage}
      \usepackage{graphicx}
      \usepackage{mathrsfs} 			
      \usepackage{bm} 				
      \usepackage{cleveref}
      \usepackage{enumitem}
      \usepackage{moreenum}			
      \usepackage{blindtext}
      \usepackage{tikz-cd}
      \usepackage{tikz}
      \usepackage{subfiles}
      \usetikzlibrary{shapes.geometric}
      \usetikzlibrary{automata, bending,matrix,arrows,decorations.pathmorphing,backgrounds,positioning,fit,petri,arrows.meta}
      \usepackage{textcomp}
      \usepackage{indentfirst}
      \usetikzlibrary{arrows}
      \usepackage{enumitem}
      \tikzset{commutative diagrams/.cd,arrow style=tikz,diagrams={>=latex'}}
      \usepackage{filecontents}
      \usepackage[alphabetic,lite]{amsrefs} 	
      \usepackage{stmaryrd}	
      \usepackage{subfiles}
      \usepackage{ragged2e}
      \usetikzlibrary{arrows}
      \usepackage{stackengine}

\usepackage[colorinlistoftodos,prependcaption,textsize=tiny]{todonotes}
\newcommandx{\unsure}[2][1=]{\todo[linecolor=red,backgroundcolor=red!25,bordercolor=red,#1]{#2}}
\newcommandx{\change}[2][1=]{\todo[linecolor=blue,backgroundcolor=blue!25,bordercolor=blue,#1]{#2}}
\newcommandx{\info}[2][1=]{\todo[linecolor=OliveGreen,backgroundcolor=OliveGreen!25,bordercolor=OliveGreen,#1]{#2}}
\newcommandx{\improvement}[2][1=]{\todo[linecolor=Plum,backgroundcolor=Plum!25,bordercolor=Plum,#1]{#2}}
\newcommandx{\thiswillnotshow}[2][1=]{\todo[disable,#1]{#2}}

\usepackage{marginnote}
\newcommand{\mytodo}[2][]{{%
 \let\marginpar\marginnote
 \reversemarginpar
 \renewcommand{\baselinestretch}{0.8}%
 \todo[#1]{#2}}}

         \newcommand{\GG}{\Gamma}
         
         \newcommand{\GGL}{\Lambda}

         \newcommand{\bA}{\mathbb{A}}

         \newcommand{\bF}{\mathbb{F}}
         \newcommand{\bG}{\mathbb{G}}

         \newcommand{\bP}{\mathbb{P}}

         \newcommand{\bZ}{\mathbb{Z}}

         \newcommand{\bbB}{\mathbf{B}}

         \newcommand{\bbZ}{\mathbf{Z}}

         \newcommand{\cC}{\mathcal{C}}
         \newcommand{\cD}{\mathcal{D}}
         \newcommand{\cE}{\mathcal{E}}
         \newcommand{\cF}{\mathcal{F}}

         \newcommand{\cO}{\mathcal{O}}
         \newcommand{\cP}{\mathcal{P}}

         \newcommand{\cS}{\mathcal{S}}

         \newcommand{\fm}{\mathfrak{m}}
         
         \newcommand{\fp}{\mathfrak{p}}

         \newcommand{\sF}{\mathscr{F}}
         \newcommand{\sG}{\mathscr{G}}

         \newcommand{\sO}{\mathscr{O}}

         \newcommand{\sX}{\mathscr{X}}
         \newcommand{\sY}{\mathscr{Y}}

         \newcommand{\ra}{\rightarrow}

         \newcommand{\hra}{\hookrightarrow}

         \newcommand{\bDelta}{\mathbf{\Delta}}
         \newcommand{\wt}{\widetilde}
         
         \newcommand{\Sing}{\mathrm{Sing}^{\bA^1}}

         \newcommand{\pr}{^{\prime}}

         \newcommand{\ce}{\colonequals}
         \newcommand{\ov}{\overline}

         \newcommand{\Fun}{\mathrm{Fun}}

         \newcommand{\isoto}{\overset{\sim}{\longrightarrow}}
         
         \newcommand{\Sm}{{\mathrm{Sm}}}  
         \newcommand{\fppf}{\mathrm{fppf}}		                                           
         \newcommand{\et}{\mathrm{\acute{e}t}}	                                           
         \newcommand{\Zar}{\mathrm{Zar}}		                                               
          \newcommand{\Shv}{\mathrm{Shv}}		                                          

\newextarrow{\xbigtoto}{{15}{15}{15}{12}}
   {\bigRelbar\bigRelbar{\bigtwoarrowsleft\rightarrow\rightarrow}}
         
         \DeclareMathOperator{\Spec}{Spec}		                       
         \DeclareMathOperator{\Spc}{Spc}
         \DeclareMathOperator{\Frac}{Frac}		                       
         \DeclareMathOperator{\Map}{Map}			                       
          
         \DeclareMathOperator{\Set}{\mathbf{Set}}
         \DeclareMathOperator{\Sch}{\mathbf{Sch}}		                                                  
         \DeclareMathOperator{\colim}{colim}		                                                  
         \DeclareMathOperator{\GL}{GL}		                                                  
         \DeclareMathOperator{\codim}{codim}		                                                  

         \newcommand{\ba}{\begin{aligned}}
         \newcommand{\ea}{\end{aligned}}
         \newcommand{\be}{\begin{equation}}
         \newcommand{\ee}{\end{equation}}
         \newcommand{\pf}{\begin{proof}}
         \newcommand{\bpf}{\begin{proof}}
         \newcommand{\epf}{\end{proof}}
         \newcommand{\bsol}{\begin{solution}}
         \newcommand{\esol}{\end{solution}}
         \newcommand{\bthm}{\begin{thm}}
         \newcommand{\ethm}{\end{thm}}
         \newcommand{\bthmt}{\begin{thm-tweak}}
         \newcommand{\ethmt}{\end{thm-tweak}}
         \newcommand{\bprop}{\begin{prop}}
         \newcommand{\eprop}{\end{prop}}
         \newcommand{\bcor}{\begin{cor}}
         \newcommand{\ecor}{\end{cor}}

         \newcommand{\brem}{\begin{rem}}
         \newcommand{\erem}{\end{rem}}
         \newcommand{\bremt}{\begin{rem-tweak}}
         \newcommand{\eremt}{\end{rem-tweak}}
         \newcommand{\brems}{\begin{rems} \hfill \begin{enumerate}[label=\b{\thenumberingbase.},ref=\thenumberingbase]}
         \newcommand{\remi}{\addtocounter{numberingbase}{1} \item}
         \newcommand{\erems}{\end{enumerate} \end{rems}}
         \newcommand{\begs}{\begin{egs} \hfill \begin{enumerate}[label=\b{\thenumberingbase.},ref=\thenumberingbase]}
         \newcommand{\egi}{\addtocounter{numberingbase}{1} \item}
         \newcommand{\eegs}{\end{enumerate} \end{egs}}
         \newcommand{\eremstweak}{\end{enumerate} \end{rems-tweak}}
         \newcommand{\eremst}{\end{enumerate} \end{rems-tweak}}
         \newcommand{\blem}{\begin{lemma}}
         \newcommand{\elem}{\end{lemma}}
         \newcommand{\blemt}{\begin{lemma-tweak}}
         \newcommand{\elemt}{\end{lemma-tweak}}
         \newcommand{\bconj}{\begin{conj}}
         \newcommand{\econj}{\end{conj}}
         \newcommand{\bprob}{\begin{Problem}}
         \newcommand{\eprob}{\end{Problem}}
         \newcommand{\bpropt}{\begin{prop-tweak}}
         \newcommand{\epropt}{\end{prop-tweak}}
         \newcommand{\bq}{\begin{Q}}
         \newcommand{\eq}{\end{Q}}
         \newcommand{\benum}{\begin{enumerate}[label={{\upshape(\alph*)}}]}
         \newcommand{\benuma}{\begin{enumerate}[label={{\upshape(\arabic*)}}]}
         \newcommand{\benumr}{\begin{enumerate}[label={{\upshape(\roman*)}}]}
         \newcommand{\eenum}{\end{enumerate}}
         \newcommand{\bc}{}
         \newcommand{\bd}{\begin{defn}}
         \newcommand{\ed}{\end{defn}}
         \newcommand{\bque}{\begin{que}}
         \newcommand{\eque}{\end{que}}
         \newcommand{\bfct}{\begin{fact}}
         \newcommand{\efct}{\end{fact}}

         \newcommand{\beg}{\begin{eg}}
         \newcommand{\eeg}{\end{eg}}

         \newcommand{\bcl}{\begin{claim}}
         \newcommand{\ecl}{\end{claim}}
         \newcommand{\bclt}{\begin{cl-tweak}}
         \newcommand{\eclt}{\end{cl-tweak}}

         \newcommand{\x}{\text}

         \newcommand{\q}{\quad}
         
         \newcommand{\qqq}{\quad\quad\quad}

         \newcommand{\tst}{\textstyle}

         \newcommand{\Nis}{\mathrm{Nis}}

\newcommand*{\QED}{\hfill\ensuremath{\qed}}

\tikzset{
    labl/.style={anchor=south, rotate=90, inner sep=.5mm}
}

      \usepackage{aliascnt}
\newaliascnt{numberingbase}{subsection}
\numberwithin{equation}{numberingbase}

\newtheoremstyle{thms}{0pt}{0pt}{\itshape}{}{\bfseries}{.}{ }{}
\theoremstyle{thms}
\newtheorem{conj}[numberingbase]{Conjecture}

\newtheorem{cor}[numberingbase]{Corollary}
\newtheorem{lemma}[numberingbase]{Lemma}

\newtheorem{prop}[numberingbase]{Proposition}

\newtheorem{Q}[numberingbase]{Question}
\newtheorem{thm}[numberingbase]{Theorem}

\newtheoremstyle{claims}{0pt}{0pt}{}{}{\itshape}{.}{ }{}
\theoremstyle{claims}
\newtheorem{claim}[equation]{Claim}
\newtheorem{cl-tweak}[subsubsection]{Claim}
\Crefname{cl-tweak}{Claim}{Claims}

\newtheoremstyle{defs}{0pt}{0pt}{}{}{\bfseries}{.}{ }{}
\theoremstyle{defs}
\newtheorem{defn}[numberingbase]{Definition}

\newtheorem{eg}[numberingbase]{Example}

\newtheorem*{egs}{Examples}
\newtheorem{rem}[numberingbase]{Remark}

\newtheorem*{rems}{Remarks}

\Crefname{claim}{Claim}{Claims}
\Crefname{sublemma}{Lemma}{Lemmas}
\Crefname{conj}{Conjecture}{Conjectures}
\Crefname{cor}{Corollary}{Corollaries}
\Crefname{defn}{Definition}{Definitions}
\Crefname{eg}{Example}{Examples}
\Crefname{prop}{Proposition}{Propositions}
\Crefname{Q}{Question}{Questions}
\Crefname{rem}{Remark}{Remarks}
\Crefname{thm}{Theorem}{Theorems}
\Crefname{variant}{Variant}{Variants}

\theoremstyle{thms}
\newtheorem{thm-tweak}[subsection]{Theorem}
\Crefname{thm-tweak}{Theorem}{Theorems}
\newtheorem{lemma-tweak}[subsection]{Lemma}
\Crefname{lemma-tweak}{Lemma}{Lemmas}
\newtheorem{cor-tweak}[subsection]{Corollary}
\Crefname{cor-tweak}{Corollary}{Corollaries}
\newtheorem{prop-tweak}[subsection]{Proposition}
\Crefname{prop-tweak}{Proposition}{Propositions}
\newtheorem{conj-tweak}[subsection]{Conjecture}
\Crefname{conj-tweak}{Conjecture}{Conjectures}

\theoremstyle{defs}
\newtheorem{defn-tweak}[subsection]{Definition}
\Crefname{defn-tweak}{Definition}{Definitions}
\newtheorem{eg-tweak}[subsection]{Example}
\Crefname{eg-tweak}{Example}{Examples}
\newtheorem*{rems-tweak}{Remarks}
\newtheorem{rem-tweak}[subsection]{Remark}
\Crefname{rem-tweak}{Remark}{Remarks}

\newtheoremstyle{subsection-tweak}
   {0pt}
   {0pt}%
   {}
   {}%
   {\bfseries}
   {}%
   {.5em}
   {\thmnumber{\@{#1}{}\@{#2}.}%
    \thmnote{~{\bfseries#3.}}}

\theoremstyle{subsection-tweak}
\newtheorem{pp}[subsection]{}
\newcommand{\bpp}{\begin{pp}}
\newcommand{\epp}{\end{pp}}
\newtheorem{pp-t}[subsubsection]{}

\theoremstyle{subsection-tweak}

\theoremstyle{subsection-tweak}
\newtheorem{pp-tweak}{}

      \makeatletter
      \def\@tocline#1#2#3#4#5#6#7{
          \begingroup
          \@ifempty{#4}{}{}

          \parindent\z@ \leftskip#3\relax \advance\leftskip\@tempdima\relax
          #5\hskip-\@tempdima
            \ifcase #1
             \or\or \hskip 2em \or \hskip 1em \else \hskip 3em \fi%
            #6\nobreak\relax
          \dotfill\hbox to\@pnumwidth{\@tocpagenum{#7}}\par
          \nobreak
          \endgroup
        }
       \def\l@section{\@tocline{1}{0pt}{1pc}{}{}}

      \renewcommand{\tocsection}[3]{%
        \indentlabel{\@ifnotempty{#2}{\makebox[1.3em][l]{%
          \ignorespaces#1 \bfseries{#2}.\hfill}}}\bfseries{#3}
          \vspace{-3.5pt}}

      \renewcommand{\tocsubsection}[3]{%
        \indentlabel{\@ifnotempty{#2}{\hspace*{-0.5em}\makebox[2.1em][l]{%
          \ignorespaces#1#2.\hfill}}}#3
          \vspace{-4.5pt}}

\makeatother

\makeatletter 
\newcommand\appendix@section[1]{%
  \refstepcounter{section}%
  \orig@section*{Appendix \@Alph\c@section. #1}%
}
\let\orig@section\section
\g@addto@macro\appendix{\let\section\appendix@section}
\makeatother

      \DeclareMathVersion{normal2}

\author{Ning Guo}
\address{Institute for Advanced Study in Mathematics of Harbin Institute of Technology, Harbin, China}
\email{ninguo@hit.edu.cn}
\author{Fei Liu}
\address{Department of Mathematics, Southern University of Science and Technology, Shenzhen, China}
\email{liufei54@pku.edu.cn}
\date{\today}

\def\UTFviii@defined#1{%
  \ifx#1\relax
      \PackageError{inputenc}{Unicode\space char\space\expandafter
                              \UTFviii@splitcsname\string#1\relax
                              \MessageBreak
                              not\space set\space up\space
                              for\space use\space with\space LaTeX}\@eha
  \else\expandafter
    #1%
  \fi
}

\makeatletter
\def\UTFviii@defined#1{%
  \ifx#1\relax
      ?%
  \else\expandafter
    #1%
  \fi
}

\makeatother


\begin{document}

\title{The Bass--Quillen conjecture for torsors over valuation rings}
\maketitle

\vspace{-25pt}

\begin{abstract}
For a valuation ring $V$, a smooth $V$-algebra $A$, and a reductive $V$-group scheme $G$ satisfying a certain natural isotropicity condition, we prove that every Nisnevich $G$-torsor on $\bA^N_A$ descends to a $G$-torsor on $A$.  
As a corollary, we generalize Raghunathan's theorem on torsors over affine spaces to a relative setting.  
We also extend several affine representability results of Asok, Hoyois, and Wendt from equi-characteristics to mixed characteristics.  
Our proof relies on previous work on the purity of reductive torsors over smooth relative curves and the Grothendieck--Serre conjecture for constant reductive group schemes.  
\end{abstract}



\hypersetup{
    linktoc=page,     
}

\renewcommand*\contentsname{}
\q\\
\tableofcontents


\section{Introduction} 

Let $A$ be a Noetherian regular ring. 
The classical Bass--Quillen conjecture, which grew out of Serre's problem that was solved by Quillen in \cite{Qui76} and Suslin in \cite{Sus76}, asserts that every vector bundle on the affine space $\bA_A^n$ descends to $A$; see \cite{Qui76}*{Comment~1 on page~170} and \cite{Bas73}*{Problem~IX}. 
Building on foundational earlier work, the conjecture was resolved affirmatively in the unramified case or, more generally, when $A$ is regular over a Dedekind ring; see \cite{Qui76}, \cite{Lin81}, and \cite{Pop89}. 
Nonetheless, the conjecture remains open in the general case.

In this article, we consider a variant of the conjecture for torsors under reductive group schemes. As will become evident, imposing a specific isotropicity condition on the reductive groups is crucial.

\bd [Total isotropicity; see \cite{Ces22a}*{Definition~8.1} or \cite{Ces22}*{Section~1.3.6}]
   Let $G$ be a reductive group scheme over a scheme $S$. 
   For a point $s\in S$, the adjoint semisimple $\sO_{S,s}$-group $G_{\sO_{S,s}}^{\text{ad}}$ decomposes as (\cite{SGA3IIInew}*{expos\'e~XXIV, proposition~5.10(i)}):
  \begin{equation}\label{decomps od G^ad}
   \tst G_{\sO_{S,s}}^{\text{ad}}\simeq \prod_{i}\text{Res}_{R_i/\sO_{S,s}}(G_i),
   \end{equation}
  where each $\Spec R_i\ra \Spec \sO_{S,s}$ is a finite \'etale cover, and $G_i$ is a simple\footnote{This means that the geometric fibres are reductive algebraic groups with connected Dynkin diagrams.} adjoint $R_i$-group scheme (of constant type). 
  We say that $G$ is \emph{totally isotropic at $s\in S$}, if each $G_i$ contains a copy of $\mathbb{G}_{m,R_i}$ as an $R_i$-subgroup. 
  We say that $G$ is \emph{totally isotropic} if it is totally isotropic at each point of $S$.
\ed

\bremt
By \cite{SGA3IIInew}*{expos\'e~XXVI, corollaire~6.12}, $G$ is totally isotropic at a point $s\in S$ if and only if every $G_i$ in (\ref{decomps od G^ad}) contains a parabolic $R_i$-subgroup that is $R_i$-fibrewise proper.
\eremt

\begs
 \egi A reductive group $G$ over a field $k$ is totally isotropic if and only if every (almost) $k$-simple factor of its derived subgroup $G^{\mathrm{der}}$ is isotropic (i.e., every factor contains a copy of $\bG_{m,k}$).
 \egi Tori and quasi-split (in particular, split) reductive group schemes over an arbitrary base scheme are totally isotropic. 
 Also, the total isotropicity is stable under base change and finite products.
\eegs

Note that Fedorov \cite{Fed21}*{Def.~1} refers to totally isotropic groups as \emph{strongly locally isotropic} groups.

The vector bundle case $G=\GL_n$ in the following conjecture recovers the original Bass--Quillen conjecture.

\bconj[Bass--Quillen for torsors] \label{Bass-Quillen for torsors}
For a Noetherian regular ring $A$ and a \emph{totally isotropic} reductive $A$-group scheme $G$, every Nisnevich-locally\footnote{Recall that a Nisnevich covering of a scheme $X$ is a family of \'etale morphisms $\{f_i\colon U_i\ra X\}_{i\in I}$ such that for every $x\in X$, there is a $u\in U_i$ lying over $x$ inducing an isomorphism of residue fields $k_x\isoto k_u$.
In particular, if $X$ is an integral scheme, every Nisnevich cover of $X$ has a generic section. If $X$ is quasi-compact and quasi-separated, then the definition above is equivalent to that the morphism $f\colon \bigsqcup_{i\in I}U_i\ra X$ admits a section over each stratum $Z_i\backslash Z_{i+1}$ for a sequence of finitely presented closed subschemes $\emptyset=Z_n\subset Z_{n-1}\subset \cdots \subset Z_1\subset Z_0=X$. 
If $X$ is Noetherian, this follows directly from a Noetherian induction argument. 
The general quasi-compact quasi-separated case is a little bit subtle and perhaps not so well known; we refer the interested reader to a note of Hoyois, see https://hoyois.app.uni-regensburg.de/papers/allagree.pdf.} trivial $G$-torsor on $\mathbb{A}_A^N$ descends to $A$, that is,
\begin{equation}\label{pullback map}
  H^1_{\mathrm{Nis}}(A,G)\xrightarrow{\sim} H^1_{\mathrm{Nis}}(\bA_A^N,G) \q \text{ via pullback.}
\end{equation}
\econj

\brems\label{remark on BQ}
\remi Both the total isotropicity of $G$ and the Nisnevich-local triviality (instead of merely \'etale-local triviality) of torsors are essential: (1) there are non-trivial \'etale $\mathrm{PGL}_n$-torsors over $\mathbb{A}_k^1$ even for a separably closed but non-algebraically closed field $k$ if $\mathrm{char}(k)|n$, so they do not descend to $k$ (see \cite{CTS21}*{théorème~5.6.1(vi)}); (2)
           there exists a generically trivial (in fact, Zariski-locally trivial) $\text{SO}(x_1^2+\cdots +x_4^2)$-torsor over the affine $\mathbf{R}$-plane $\mathbb{A}_{\mathbf{R}}^2$ that does not descend to $\mathbf{R}$ (see \cite{Par78}). In fact, Balwe and Sawant have shown in \cite{BS17}*{Theorem~4.8} that if $G$ is defined over an infinite perfect field $k$ and is not totally isotropic, then the pullback map (\ref{pullback map}) cannot be bijective for all smooth $k$-algebras. See also \cite{Fed16} for more concrete counterexamples.
\remi \label{rem on Nis vs Zar} The Grothendieck--Serre conjecture predicts that an \'etale torsor under a reductive group scheme $G$ over a regular semilocal domain $A$ is trivial if it trivializes over the fraction field. 
This conjecture was settled affirmatively when $A$ contains a field (\cite{FP15} and \cite{Pan20a}), when $A$ is a semilocal Dedekind domain $R$ (see \cite{Nis82, Nis84} and \cite{Guo22}; this case, by induction, implies the case when $A$ is Henselian).
For $A$ that is essentially smooth over $R$ (namely, $A$ is \emph{unramified}), the subcases when $G$ descends to a reductive $R$-group and when $G$ is totally isotropic are resolved in \cite{GL23}*{\S~1.3} and \cite{CF23} respectively. 
For a comprehensive summary of the state of the art, see \cite{GL23} or \cite{CF23}. 
The last aforementioned case implies that, when $A$ is unramified in the context of \Cref{Bass-Quillen for torsors}, a $G$-torsor over $\bA_A^N$ is Nisnevich-locally trivial, if and only if it is generically trivial, if and only if it is Zariski-locally trivial. Therefore, Conjecture \ref{Bass-Quillen for torsors} can be reformulated as:
        \[
H_{\mathrm{Zar}}^1(A,G) \xrightarrow{\sim} H_{\mathrm{Zar}}^1(\mathbb{A}_{A}^N,G) \q \text{ via pullback}.
        \]
\erems

Apart from the classical vector bundle case $G=\GL_n$, \Cref{Bass-Quillen for torsors} was proved in the following cases:
\begin{itemize}
  \item   Asok, Hoyois, and Wendt settled the case when $A$ is smooth over a field $k$ and $G$ is the pullback of a totally isotropic reductive $k$-group; see \cite{AHW18} for $k$ infinite and \cite{AHW20} for $k$ finite.
  \item Stavrova established the case when $A$ contains a field in \cite{Sta22}*{Corollary~5.5}, and an earlier result \cite{Sta19}*{Theorem~4.4} addressed the case when $A$ contains an infinite field. (For the convenience of readers, a greatly simplified version of her proof is outlined in \S~\Cref{pf-of-equi-char-BQ-torsors}.)
  \item  $\check{\mathrm{C}}$esnavi$\check{\mathrm{c}}$ius, in \cite{Ces22b}, independently proved \Cref{Bass-Quillen for torsors} in the case when $A$ contains a field.
  This was achieved by completing the solution of Nisnevich's purity conjecture over such an $A$, as provided by Fedorov \cite{Fed21}, and subsequently deriving it from Nisnevich's purity conjecture.
\end{itemize}

Consequently, the equi-characteristic case of \Cref{Bass-Quillen for torsors} has been fully settled. 
However, to the best of our knowledge, the mixed characteristic case remains open.
Moreover, as readers will note in \S~\Cref{pf-of-equi-char-BQ-torsors}, the approach used for the equi-characteristic case does not readily extend to the mixed characteristic case at present, primarily because the relative Grothendieck--Serre conjecture (\Cref{mixed-char-rel-G-S}) remains unresolved.

The first main result of this article settles a variant of \Cref{Bass-Quillen for torsors}.
This result is novel even in the specific case where $R$ is a discrete valuation ring (DVR) of mixed characteristic, stated as follows.
\bthmt [\Cref{B-Q over val rings}] \label{intro-B-Q over val rings-dvr}
  Let $A$ be a ring that is ind-smooth over a discrete valuation ring $R$, and let $G$ be a totally isotropic reductive $R$-group scheme. Then, via pullback, we have the following bijection
\[
\text{$H^1_{\mathrm{Nis}}(A,G)\isoto H^1_{\mathrm{Nis}}(\bA^N_A,G),$\q or, equivalently, \q $H^1_{\mathrm{Zar}}(A,G)\isoto H^1_{\mathrm{Zar}}(\bA^N_A,G)$}.
\]
\ethmt

Moreover, we investigate \Cref{Bass-Quillen for torsors} for a broad class of rings $A$ that resemble Noetherian regular rings, specifically those that are (ind-)smooth over a Prüfer ring, and thus may be non-Noetherian. 
Recall that a ring is \emph{Prüfer} if all its local rings are valuation rings; a \emph{valuation ring} is a domain whose ideals are totally ordered by inclusion.
While Noetherian Prüfer rings correspond precisely to Dedekind rings, non-Noetherian Prüfer rings are equally prevalent: they arise naturally in non-Archimedean analytic geometry and perfectoid theory, with the integer rings of perfectoid fields serving as notable examples.

More significantly, the recently introduced v- and arc-topologies define covers by testing maps from (spectra of) valuation rings, providing a powerful framework for reducing various problems to analogous ones for valuation rings. 
Furthermore, growing evidence suggests that a ring smooth over a Pr\"ufer ring exhibits behaviour akin to that of a Noetherian regular local ring; see \cite{GL24} for further discussion of this perspective.
These insights motivate the study of geometry over Pr\"ufer rings in \emph{loc.~cit.} and, in our scenario, the Bass--Quillen \Cref{Bass-Quillen for torsors} within this broader framework.

\bthmt [\Cref{B-Q over val rings}] \label{intro-B-Q over val rings}
  Let $A$ be a ring that is ind-smooth over a Pr\"{u}fer ring $R$, and let $G$ be a totally isotropic reductive $R$-group scheme. Then, via pullback, we have the following bijection
\[
\text{$H^1_{\mathrm{Nis}}(A,G)\isoto H^1_{\mathrm{Nis}}(\bA^N_A,G),$\q or, equivalently, \q $H^1_{\mathrm{Zar}}(A,G)\isoto H^1_{\mathrm{Zar}}(\bA^N_A,G)$}.
\]
\ethmt

\bremt
 In view of this result, we expect that \Cref{Bass-Quillen for torsors} holds for all totally isotropic reductive group scheme over a ring $A$ (instead of merely over $R$) that is ind-smooth over a  Pr\"ufer ring.
\eremt

In the case $G=\GL_n$, we note that special subcases of \Cref{intro-B-Q over val rings} have already been studied in the literature. Specifically, Simis and Vasconcelos considered the case when $A$ is a valuation ring and $N=1$ in \cite{SV71}, while Lequain and Simis treated the case when $A$ is a Pr\"{u}fer ring in \cite{LS80} (its special case when $A$ is a DVR was previously settled in \cite{Qui76}*{Theorem~4'}). 
Insofar as we are aware, there are no other instances of \Cref{intro-B-Q over val rings} in such a non-Noetherian context, even for the simplest group $G=\GL_n$.

The proof of \Cref{intro-B-Q over val rings} follows the same spirit as the approach to the Bass--Quillen conjecture for vector bundles in \cite{Lin81} and \cite{Pop89}, but the technical details are significantly more intricate and demanding. 
This proof is enabled by recent progress on the mixed characteristic Grothendieck--Serre conjecture obtained in \cite{GL23}, particularly the new geometric presentation lemma (\Cref{variant of Lindel's lem}), and the Pr\"uferian purity result for $G$-torsors (\Cref{purity for rel. dim 1}), specifically its important corollary, \Cref{triviality over R(t)}.

 \bpp[Torsors over relative affine spaces]
 In \cite{Rag89}*{Theorem~A}, Raghunathan proved that for a totally isotropic reductive group $G$ over a field $k$, a $G$-torsor over $\bA^N_k$ is trivial if it is trivial over $\bA^N_{k^s}$ and over the origin $0\in \bA^N_k(k)$, where $k^s$ denotes a separable closure of $k$. As an application of \Cref{intro-B-Q over val rings}, we can efficiently reprove this result and even present a generalization, as stated below.
\epp

\bthmt\label{Raghunathan-introd}
Let $A$ be a domain with fraction field $K$ and $G$ a reductive $A$-group scheme.
Assume 
\begin{itemize}
  \item either $A$ is regular Noetherian and contains a field, and $G$ is totally isotropic,
  \item or $A$ is ind-smooth over a Pr\"{u}fer ring $ R $, and $ G$ descends to a totally isotropic reductive $R$-group.
\end{itemize}
  Then, a $G$-torsor over $\bA^N_A$ is trivial if and only if it is trivial over $\bA^N_{K^s}$ and the zero section $0_A\in \bA^N_A(A)$.
\ethmt

\bremt
Assuming the most general case of \Cref{Bass-Quillen for torsors}, our argument will extend \Cref{Raghunathan-introd} to the case of an arbitrary ind-smooth domain $A$ over a Pr\"{u}fer domain. 
\eremt

\bpp[$\bA^1$-homotopy theoretic affine representability in mixed characteristics]
As applications of \Cref{intro-B-Q over val rings}, we can generalize many results from \cite{AHW18} and \cite{AHW20}, extending them beyond equi-characteristics to encompass mixed characteristics as well.

Let $S$ be a quasi-compact quasi-separated scheme. Denote by $\Spc(S)$ the $\infty$-category of motivic spaces over $S$ (it is overviewed in \S~\ref{overview-Motivic-spaces}). 
If $\sX$ and $\sY$ are presheaves of spaces on the category of $S$-smooth schemes, then we write 
$$
[\sX,\sY]_{\bA^1}\ce \pi_0 \Map_{\Spc(S)}(L_{\mathrm{mot}}(\sX), L_{\mathrm{mot}}(\sY)),
$$ 
where $L_{\mathrm{mot}}\colon \cP(\Sm_S)\to \Spc(S)$ denotes the motivic localization functor (see \S~\ref{overview-Motivic-spaces}).  Given an $S$-group scheme $G$, write $\bbB G$ for its usual bar construction\footnote{Specifically, $\bbB G$ is the functor that sends an $S$-smooth scheme $U$ to the nerve of $G(U)$, where $G(U)$ is treated as a category with a single object whose endomorphisms are given by $G(U)$.}, considered as a presheaf of spaces on the category of smooth $S$-schemes. 
We establish the following representability of Nisnevich-locally trivial torsors.
\epp

\bthmt \label{affine-representability of H^1Nis}
Let $S$ be the spectrum of a Pr\"ufer ring (e.g., a DVR) and let $G$ be a totally isotropic reductive $S$-group scheme. 
Then, for every smooth affine $S$-scheme $U$, there is a functorial bijection 
\[
H^1_{\mathrm{Nis}}(U,G) \simeq [U,\bbB G]_{\bA^1}.
\]
\ethmt

In \cite{AHW18}*{Theorem~4.1.3} and \cite{AHW20}*{Theorem~2.5}, \Cref{affine-representability of H^1Nis} was established under the assumption that \(S\) is the spectrum of a field. 
In fact, using \cite{AHW18}*{Theorem~2.3.5}, \Cref{affine-representability of H^1Nis} can be deduced from the assertion that \(\bbB_{\Nis}G\) is \(\bA^1\)-naive, as defined in \cite{AHW18}*{Definition~2.1.1} (see \S~\ref{Represn-results}). 
This assertion, in turn, is a direct consequence of the Bass--Quillen \Cref{intro-B-Q over val rings}.

As an application of \Cref{affine-representability of H^1Nis}, we obtain the following result, which, for certain homogeneous spaces of totally isotropic reductive groups, identifies naive homotopy classes and
 true $\bA^1$-homotopy classes.

\bthmt \label{naive-vs-A1-homotopy-G/H}
Let $S$ be the spectrum of a Pr\"ufer ring (e.g., a DVR). 
Suppose that
\benumr 
\item $H\to G$ is a closed immersion of totally isotropic reductive $S$-group schemes, and
\item the \'etale $H$-torsor $G\to G/H$ is Nisnevich-locally trivial\footnote{Equivalently, generically trivial by \cite{GL23}*{Theorem 1.1}.}.
\eenum
Then, $G/H$ is $\bA^1$-naive. 
In particular, for every smooth affine $S$-scheme $U$, there is a functorial bijection 
\[
\pi_0 \left( \Sing G/H\right)(U) \simeq [U, G/H]_{\bA^1}.
\]
\ethmt
 \Cref{naive-vs-A1-homotopy-G/H}  was previously established in \cite{AHW18} and \cite{AHW20} under the assumption that $S$ is the spectrum of a field.



\bpp[Notations and conventions]
Throughout this article, we work with commutative rings with units.
We freely use the language of $\infty$-categories, as set out in \cite{HTT} and \cite{HA}. By abuse of notation, the nerve of an ordinary category $\cC$ will also be denoted by $\cC$, unless otherwise specified. 
\begin{itemize}
  \item For a point $s$ of a scheme (resp., for a prime ideal $\mathfrak{p}$ of a ring), $k_s$ (resp., $k_{\mathfrak{p}}$) denotes its residue field. 
  The total ring of fractions of a ring $A$ is denoted by $\Frac A$.
      The base change along a morphism of schemes $S\ra S\pr$ is denoted by $(-)_{S\pr}$; if $S\pr=\Spec R\pr$ is affine, we also write $(-)_{R\pr }$.
   \item For a quasi-compact quasi-separated scheme $S$, we let $\mathrm{Sm}_S\ce \mathrm{Sm}_S^{\mathrm{qcqs}}$ denote the category of quasi-compact quasi-separated smooth $S$-schemes.
  \item  Let $S$ be a scheme, and let $G$ be an $S$-group scheme. Given an $S$-scheme $T$, a $G$-torsor over $T$ refers to a $G_T\ce G\times_ST$-torsor. If $G$ is $S$-smooth (typically reductive in this article) and $\tau \in \{\Nis, \et \}$,  $\bbB_{\tau} G$ denotes the classifying stack of $G$-torsors on the (small, big, etc.) $\tau$-site of $S$.
  \item $\bDelta$ denotes the simplex category, i.e., the category of finite non-empty linearly ordered sets.
  \item $\cS$ denotes the $\infty$-category of spaces (or anima, as per \cite{CS24}*{\S~5.1}).
  \item For $\infty$-categories $\cC$ and $\cD$, let $\cP(\cC, \cD)\ce \Fun(\cC^{\mathrm{op}},\cD)$ denote the $\infty$-category of $\cD$-valued presheaves on $\cC$. For brevity, we write $\cP(\cC)$ for $\cP(\cC,\cS)$. If $\tau$ is a Grothendieck topology on $\cC$, we let $\Shv_{\tau}(\cC) \subset \cP(\cC)$ denote the full subcategory of $\cS$-valued $\tau$-sheaves.
  \item $\Delta^{\bullet}_{\bbZ}$ denotes the standard cosimplicial scheme over $\Spec \bbZ$, i.e., the functor $\bDelta\ra \mathbf{Sch}_{/\mathbf{Z}}$  given by
      \[
      [n]\mapsto \Delta^n_{\bbZ}\ce \Spec \left(  \frac{\bbZ[x_0,\cdots, x_n]}{(\sum_{i=0}^nx_i-1)}   \right).
      \]
      For any scheme $S$, we let $\Delta^{\bullet}_S\colon \bDelta\ra \mathbf{Sch}_{/S}$ denote the base change of $\Delta^{\bullet}_{\bbZ}$ to $S$.
  \item $\Sing$ denotes the singular construction endofunctor on $\cP(\mathrm{Sm}_S)$, see (\ref{Singular-constru}).
\end{itemize}

\epp

\subsection*{Acknowledgements}
We would like to express our gratitude to K\k{e}stutis \v{C}esnavi\v{c}ius for his constant encouragement throughout this project.
Also, his insightful comments pointed out the inaccuracies of the article for us to improve it.
This project benefited from the first author's privileged working environment at the Institute of Advanced Sciences for Mathematics, Harbin Institute of Technology.

\section{Purity for torsors under reductive groups}
\label{sect-purity of reductive torsors}

The primary aim of this section is to present the purity of reductive torsors on smooth relative curves over Pr\"ufer bases, as stated in \Cref{purity for rel. dim 1}. 
We then derive the main result, \Cref{triviality over R(t)}.
Finally, we recall the statement of Grothendieck--Serre conjecture for constant reductive group schemes in \Cref{G-S for constant reductive gps}.

To help the reader get a feeling for working over Pr\"ufer bases, we recall the following ring-theoretic result.

\blemt \label{geom}
 Let $X$ be a scheme that is flat and locally of finite type over an integral Pr\"ufer scheme $S$.
\benumr
  \item \label{coherence of X}The scheme $X$ is locally of finite presentation and locally coherent.
  \item \label{coherence of O}For every point $x\in X$, the local ring $\sO_{X,x}$ is coherent.
  \item\label{geo-i} If $X$ is irreducible, then all nonempty $S$-fibres have the same dimension.
  \item\label{geo-iii} If $\sO_{X_s,\xi}$ is reduced for a maximal point $\xi\in X_s$, then the local ring $\sO_{X,\xi}$ is a valuation ring and the extension $\sO_{S,s}\hra \sO_{X,\xi}$ induces an isomorphism of value groups.
 \eenum
\elemt
\bpf
For \ref{coherence of X}--\ref{coherence of O}, see \cite{GL24}*{Lemma~3.2.1}. For \ref{geo-i}, see \cite{EGAIV3}*{lemme~14.3.10}.
For \ref{geo-iii}, see \cite{MB22}*{théorème~A}.
\epf

The following result, combined with limit arguments,
often allows one to only consider Pr\"ufer rings of finite Krull dimension.
\blemt \label{approxm semi-local Prufer ring}
Every semilocal Pr\"{u}fer domain $R$ is a filtered direct union of its subrings $ R_i$ such that:
\benumr
\item for every $i$, the ring $R_i$ is a semilocal Pr\"{u}fer domain of finite Krull dimension; and
\item for $i$ large enough, $R_i\to R$ induces a bijection on the sets of maximal ideals, hence is fpqc.
\eenum
\elemt
\bpf
Write $\Frac (R)=\cup_i K_i$ as the filtered direct union of the subfields of $\Frac (R)$ that are finitely generated over its prime field $\mathfrak{K}$.
For $R_i\ce R\cap K_i$, we have $R=\cup_i R_i$.
For (i), it suffices to see that every $R_i$ is a semilocal Pr\"{u}fer domain whose local rings have finite ranks.
Let $\{\fp_j\}_{1\le j \le n}$ be the set of maximal ideals of $R$. Then $R=\bigcap_{1\le j \le n} R_{\mathfrak{p}_j}$ is the intersection of the valuation rings $R_{\mathfrak{p}_j}$. Thus we have
\[
   \tst R_i=\bigcap_{1\le j \le n} \left(K_i \cap R_{\mathfrak{p}_j}\right).
\]
 Since $K_i/\mathfrak{K}$ has finite transcendence degree, by Abhyankar's inequality \cite{Abh56}*{Corollary~1}, every $K_i \cap R_{\mathfrak{p}_j}$ is a valuation ring of finite rank. By
 \cite{Bou98}*{VI, \S~7, proposition~1--2}, $R_i$ is a semilocal Pr\"{u}fer domain, and its local rings at maximal ideals are precisely the minimal elements of the set $\{K_i \cap R_{\mathfrak{p}_j}\}_{1\le j \le n}$ under inclusion.  
 This implies (i). 
 For (ii), the quasi-compact property follows from the affineness of $\Spec R\ra \Spec R_i$; each map $R_i\ra R$ is flat since $R$ is torsion free.
 As flat morphism lifts generalizations, it remains to show the bijection between maximal ideals of $R_i$ and of $R$. 
Namely, we show that for $i$ large enough there is no strict inclusion relation between $K_i \cap R_{\mathfrak{p}_{j_1}}$ and $K_i \cap R_{\mathfrak{p}_{j_2}}$ for $j_1\neq j_2$. Indeed, if $\pi_j \in \mathfrak{p}_j\backslash \bigcup_{j'\neq j} \mathfrak{p}_{j'}$ for $1\le j \le n$, then (ii) holds for any $i$ for which $\{\pi_j\}_{1\le j \le n} \subset K_i$.
 \epf

The following result will facilitate the passage to the case of closed point in the closed fibre.
\blem\label{enlarge valuation rings}
 Let $V$ be a valuation ring and let $V\ra A$ be an essentially finitely presented (resp., essentially smooth) local homomorphism of local rings. 
 There are an extension of valuation rings $V\subset V\pr$ identifying their value groups and an essentially finitely presented (resp., essentially smooth) homomorphism $V\pr \to A$ of $V$-algebras inducing a finite residue field extension.
\elem
\bpf
Assume that $A=\sO_{X,x}$ for an affine scheme $X$ finitely presented over $V$ and a point $x\in X$ lying over the closed point $s \in \Spec(V)$. Let $t=\mathrm{tr.deg}(k_x/k_s)$.
After shrinking $X$ around $x$ if necessary, choose sections $b_1,\ldots,b_t\in \GG(X,\sO_X)$ such that their images $\ov{b_1},\cdots, \ov{b_t}$ in $k_x$ form a transcendental basis of $k_x/k_s$.
Define $p\colon X\to \mathbb{A}_V^t$ by sending the standard coordinates $T_1,\ldots, T_t$ of $\mathbb{A}_V^t$ to $b_1,\ldots, b_t$, respectively.
Since the transcendental degree of $k_s(\ov{b_1},\cdots, \ov{b_t})/k_s$ is $t$, the image $\eta\ce p(x)$ is the generic point of $\bA^t_{k_s}$, so $V\pr\ce \sO_{\bA^t_{V},\eta}$ is a valuation ring whose value group is $\GG_{V\pr}\simeq \GG_V$.
Note that $k_x/k_{\eta}$ is finite, the map $V\pr\ra A$ induces a finite residue field extension.

Now assume that $V\ra A$ is essentially smooth.
For a closed point $z\in X$ specializing $x$, after shrinking $X$ around $z$, there is a factorization $X\overset{\pi}{\ra} \bA^N_V\ra \Spec V$ for an \'etale morphism $\pi$.
As the \'etale morphism $\pi$ induces finite separable extension of residue fields, we may assume that $V=k_s$ and $X=\bA^N_s$.
Since $\ov{\{x\}}$ is an irreducible closed subscheme of $X=\bA^N_s$, either we are done or it suffices to take the projection $\bA^N_s\ra \bA^{N-1}_s$ given by a standard coordinate iterately until the image of $x$ is the generic point of $\bA^t_s$. 
\epf

The following result provides an analog of the purity theorem of Colliot-Th\'el\`ene and Sansuc for reductive group torsors over two-dimensional regular Noetherian schemes (\cite{CTS79}*{th\'eor\`eme 6.13}).
\bthmt[\cite{GL24}*{Theorem~6.3}]\label{purity for rel. dim 1}
Let $X$ be a smooth relative curve over a semilocal affine Pr\"{u}fer scheme $S$. Let $Z\subset X$ be a closed subscheme such that the inclusion $j\colon X\backslash Z\hra X$ is quasi-compact, and
\[
\text{$Z_{\eta}= \emptyset$\q for each generic point $\eta\in S$ \q and \q  $\codim(Z_s,X_s)\ge 1$ \q for all $s\in S$.}
\]
Let $G$ be a reductive $X$-group scheme. 
Then, restriction induces the following equivalence 
\begin{equation*}\label{restriction of G-torsors}
  \bbB_{\et}G(X) \isoto \bbB_{\et}G(X\backslash Z).
\end{equation*}
In particular, we have a bijection of pointed sets
$$
H^1_{\et}(X,G)\simeq H^1_{\et}(X\backslash Z,G).
$$
\ethmt

\brem \label{rem on purity fo rel curves}
In higher relative dimensions, the purity \Cref{purity for rel. dim 1} is inapplicable, even in the Noetherian setting and for $G=\GL_{n}$. 
For instance, for any Noetherian regular local ring $(R,\mathfrak{m}_R)$ of Krull dimension at least $3$, there exists a vector bundle over $\Spec(R)\backslash \{\mathfrak{m}_R\}$ that cannot be extended to $\Spec R$.
\erem

The following resolution of the Grothendieck--Serre conjecture for constant reductive group schemes is a key input for our \Cref{intro-B-Q over val rings}.
In particular, we first use it to deduce \Cref{triviality over R(t)}.

\bthmt [\cite{GL23}*{Theorem~1.1}] 
\label{G-S for constant reductive gps}
For a Pr\"ufer ring $R$, an irreducible, affine, smooth $R$-scheme $X$, and a reductive $R$-group scheme $G$, every generically trivial $G$-torsor on $X$ is Zariski-semilocally trivial:
\[
   \text{the sequence \q $1\ra H^1_{\Zar}(X,G)\ra H^1_{\et}(X,G)\ra H^1_{\et}(K(X),G)$\q is exact,}
\]
where $K(X)$ is the function field of $X$.
In other words, for every semilocal ring $A$ of $X$, we have
\[
  \ker\,(H^1_{\et}(A,G)\ra H^1_{\et}(\Frac A, G))=\{\ast\}.
\]
\ethmt

For any commutative unital ring $A$, we let $A(t)$ denote the localization $S^{-1}A[t]$ with respect to the multiplicative system $S$ of \emph{monic} polynomials in $A[t]$.
\bthmt\label{triviality over R(t)}
Let $R$ be a semilocal Pr\"ufer domain and let $G$ be a reductive $R$-group scheme.
Then, a $G$-torsor over $R(t)$ is trivial if and only if it is generically trivial.
\ethmt

\bpf
Let $\cE$ be a generically trivial $G$-torsor over $R(t)$. 
Denote by $\mathfrak{r}$ the Jacobson radical of $R$. 
We observe that $R(t)$ is the semilocalization of the projective $t$-line $\mathbb{P}_R^1$ over $R$ along the infinity section $\infty_{R/\mathfrak{r}} \in \mathbb{P}_{R/\mathfrak{r}}^1$, with $s\ce \frac{1}{t}$ inverted:
\[
\tst R(t)=\cO_{\bP_R^1,\infty_{R/\mathfrak{r}}}\left[ \frac{1}{s}\right] .
\]
Note that $\infty_{R/\mathfrak{r}}$ is precisely the set of closed points of $\{s=0\}$. Hence, $\cE$ spreads out to a generically trivial $G$-torsor $\widetilde{\cE}$ on a punctured neighborhood of $\{s=0\}$ in $\bP^1_R$.
By patching torsors, we may first extend $\widetilde{\cE}$ so that its definition locus contains the generic fibre $\bP^1_{K}$, where $K=\text{Frac}(R)$. Similarly, we may assume that $\widetilde{\cE}$ is defined at the maximal point $\xi$ of each $R$-fibre of $\bP_R^1$; note that $\cO_{\bP_R^1,\xi}$ is a valuation ring. To see this, we may focus on the case where $\Spec R$ (and so $\bP_R^1$) is topological Noetherian (using \Cref{approxm semi-local Prufer ring}); then we again patch torsors using Noetherian induction and the Grothendieck--Serre conjecture for valuation rings (\cite{Guo24}*{Theorem~1.3}): every generically trivial $G$-torsor defined on an open subset of $\Spec \cO_{\bP^1_{R},\xi}$ is trivial.
Now, the definition locus of $\widetilde{\cE}$ satisfies the conditions of $X\backslash Z$ in the purity \Cref{purity for rel. dim 1}, so $\cE$ extends to a $G$-torsor $\widetilde{\cE}$ on $\bP^1_R$.
Hence, by \Cref{G-S for constant reductive gps}, the pullback of $\widetilde{\cE}$ to the semilocal ring  $\cO_{\bP_R^1,\infty_{R/\mathfrak{r}}}$ is trivial, so is the further pullback $\widetilde{\cE}|_{R(t)}=\cE$.
\epf

\section{Torsors over relative affine spaces}

\bpp[The Quillen patching and its inverse]
A central technique for studying torsors over $\bA^N_R$ is a local-to-global principle known as Quillen patching. A key insight of Gabber allows one to generalize it to rather general classes of group-valued functors, see \cite{Ces22b}*{Corollary~5.1.5}. We record it below only for locally finitely presented group schemes, which is the main case of interest for us. 
\epp

\bthmt
\label{Quillen-patching}
Let $R$ be a ring and let $G$ be a locally finitely presented $R$-group scheme.
\benum
\item \label{QP-a}
For a $G$-torsor $X$ over $R[t_1, \ldots, t_d]$, 
the set $S \subset R$ of those $r \in R$ such that $X|_{(R[t_1,\ldots, t_d])[\frac{1}{r}]}$ descends to a $G$-torsor over $R[\frac{1}{r}]$ is an ideal. 
 
\item \label{QP-b}
A $G$-torsor over $R[t_1, \ldots, t_d]$ descends to a $G$-torsor over $R$ iff it does so Zariski-locally on~$R$. 
 
 \eenum
 More generally, the analogues of \ref{QP-a} and \ref{QP-b} hold with $R[t_1, \ldots, t_d]$ replaced by any $\bZ^{\oplus d}_{\ge 0}$-graded $R$-algebra $A \cong \bigoplus_{i_1, \ldots, i_d \ge 0} A_{i_1, \ldots, i_d}$ such that $R \xrightarrow{\sim} A_{0, \ldots, 0}$. 
\ethmt

Compared to Quillen patching  \Cref{Quillen-patching}, the following `inverse' patching construction is more elementary, but still quite useful. The case where $G=\GL_n$ and $A=R[t_1,\ldots,t_N]$ is due to Roitman \cite{Roi79}*{Proposition~2}.

\blemt\label{inverse patching}
Let $R$ be a ring, let $G$ be a quasi-affine, flat, finitely presented $R$-group scheme, let $A=\bigoplus_{i_1,\ldots,i_N\ge 0}A_{i_1,\ldots,i_N}$ be a $\mathbb{Z}_{\ge 0}^{\oplus N}$-graded $R$-algebra (resp., a $\mathbb{Z}_{\ge 0}^{\oplus N}$-graded domain over $R$) such that $R\xrightarrow{\sim}A_{0,\ldots,0}$, and suppose that every $G$-torsor on $A$ (resp., every generically trivial $G$-torsor on $A$) descends to a $G$-torsor on $R$. Then, for any multiplicative subset $S\subset R$, every $G$-torsor on $A_S$ (resp., every generically trivial $G$-torsor on $A_S$) whose restriction to each local ring of $(A_{0,\ldots,0})_S\simeq R_S$ extends to a $G$-torsor on $R$ descends to a $G$-torsor on $R_S$.
\elemt
(The relevant case for us is when $A=R[t_1,\ldots,t_N]$.)

\bpf
  We focus on the part on generically trivial torsors, since the other is \cite[Proposition~5.1.10]{Ces22}.

  Let $X$ be a generically trivial $G$-torsor on $A_S$ whose restriction to each local ring of $(A_{0.\ldots,0})_S\simeq R_S$ extends to a $G$-torsor on $R$. Using Quillen patching \Cref{Quillen-patching}, we can enlarge $S$ to reduce to the case when $R_S$ is local. Then, by our assumption, the restriction of $X$ to $(A_{0.\ldots,0})_S\simeq R_S$ extends to a $G$-torsor $X_0$ on $R$. Applying a limit argument, we can reduce to the case when $S={r}$ is a singleton, at the cost of $R_S$ no longer being local. Notice that the projection onto the $(0,\ldots,0)$-th component
  \[
  \tst R\oplus \bigl( \bigoplus_{(i_1,\ldots,i_N)\neq (0,\ldots, 0)} A_{i_1,\ldots,i_N}[\frac{1}{r}]  \bigr)\simeq A[\frac{1}{r}] \times_{R[\frac{1}{r}]}R \twoheadrightarrow R
  \]
  induces an isomorphism both modulo $r^n$ and on $r^n$-torsion for every $n>0$. So, by \cite[Proposition~4.2.2]{Ces22}, we can glue the $G$-torsor $X$ on $A[\frac{1}{r}]$ with the $G$-torsor $X_0$ on $R$ to obtain a generically trivial $G$-torsor $\widetilde{X}$ on $A[\frac{1}{r}] \times_{R[\frac{1}{r}]}R$. Observe that
\[
\tst A[\frac{1}{r}] \times_{R[\frac{1}{r}]}R \simeq \underset{i\in \mathbb{N}}{\mathrm{colim}} \,A,
\]
where the transition maps $A\to A$ are maps of $R$-algebras that become isomorphisms over $R_S$ and are given by multiplications by $r^{i_1+\ldots+i_N}$ on the degree $(i_1,\ldots,i_N)$-part $A_{i_1,\ldots,i_N}$. Hence, by a standard limit argument, $\widetilde{X}$ descends to a generically trivial $G$-torsor on some copy of $A$ in the direct colimit. Therefore, by assumption, it descends further to a $G$-torsor on $R$. The base change to $R_S$ of this final descent gives a desired descent of $X$ to a $G$-torsor on $R_S$.
\epf

\bpp[Torsors on $\bA_R^N$ under reductive $R$-group schemes]
The following result was conjectured in \cite{Ces22}*{Conjecture~3.5.1} and settled later in \cite{Ces22b}*{Theorem~2.1(a)}.
\epp
\bthmt\label{triviality over relative affine space for totally isotropic}
  For a ring $R$ and a $\mathrm{totally}$ $\mathrm{isotropic}$ reductive $R$-group scheme $G$, any $G$-torsor on $\mathbb{A}_R^N$ that is trivial away from some $R$-finite closed subscheme of $\mathbb{A}_R^N$ is trivial.
\ethmt

Notice that the isotropicity assumption on $G$ is essential, see, e.g., \cite{Fed16} for counterexamples.

\begin{proof}[Sketch of proof]
Let $\cP$ be a $G$-torsor on $\mathbb{A}_R^N$ that trivializes over $\mathbb{A}_R^N\backslash Z$, where $Z\subset \mathbb{A}_R^N$ is a $R$-finite closed subscheme.
  As $Z$ is also finite under the projection $\bA^N_R\to \mathbb{A}_R^{N-1}$ onto the first $(N-1)$-coordinates, replacing $R$ with $\mathbb{A}_R^{N-1}$ reduces us to the key case when $N=1$. In this case, one extends $\cP$ to a $G$-torsor $\widetilde{\cP}$ over $\bP_R^1$ (by gluing $\cP$ with the trivial torsor over $\bP_R^1\backslash Z$). 
  Then $\widetilde{\cP}|_{\infty_R}$ is trivial, and, by base change along $\bA_R^1\to \Spec R$ (as observed by Gabber), one can reduce to showing that $s^*(\widetilde{\cP})$ is trivial for a section $s\in \bA_R^1(R)$ (in fact, even $s=0_R$ by shifting). 
  This statement is insensitive to replacing $\widetilde{\cP}$ by its pullback along the map $\bP_R^1\to \bP_R^1$ given by $[x \colon y]\mapsto [x^d\colon y^d]$ for an integer $d>0$, which allows one to assume that $G$ is moreover semisimple simply connected. 
  Granted these reductions, we revert to prove the stronger statement that $\widetilde{\cP}|_{\bA_R^1}$ is trivial. For this, Quillen patching \Cref{Quillen-patching} reduces us to the case of a local $R$. In this local case, one can try to modify $\widetilde{\cP}$ along $\infty_R$ so that it becomes trivial on the closed fibre $\bP_{R/\fm_R}^1$, and then argue using deformation rigidity of the trivial torsor on $\bP^1$.
\end{proof}

For tori, the following result shows that \Cref{Bass-Quillen for torsors} holds more generally for fppf-torsors (not merely Nisnevich-torsors) over any integral, normal base ring.
\blemt \label{B-Q for mult type}
For a normal domain $A$, an $A$-group $M$ of multiplicative type, the pullback
\[
    \text{$H^1_{\fppf}(A,M)\isoto H^1_{\fppf}(\mathbb{A}_A^n,M)$ \q is bijective. }
\]
\elemt
\bpf
When $A$ is Noetherian, this is \cite{CTS87}*{Lemma~2.4}. For a general normal domain $A$, we write it as a filtered union of its finitely generated $\mathbf{Z}$-subalgebras $A_i$, and, by replacing $A_i$ with its normalization (which is again of finite type over $\mathbf{Z}$), we may assume that each $A_i$ is normal. Then we can conclude from the Noetherian case via a limit argument, since $M$ is finitely presented over $A$.
\epf

\bpp[A Geometric presentation lemma]\label{geom_pres}
   A useful lemma due to Lindel \cite{Lin81}*{Proposition 1 and Lemma} states that an \'etale extension of local rings $A\to B$ with trivial extension of residue fields induces isomorphisms along a well-chosen nonunit $r\in A$:
   $$
   A/r^nA\isoto B/r^nB, \q \text{ where } \q n\ge1.
   $$
   In our context, where the prescribed $B$ is essentially smooth over a valuation ring, we proved a variant of \emph{loc.~cit.} that allows us to fix the $r\in B$ in advance, at the cost of the carefully choosing $A$ to be a local ring of an affine space over that valuation ring. This result is one of the crucial geometric tools for attacking the Grothendieck--Serre conjecture for `constant' reductive group schemes in \cite{GL23}. Like Lindel's work on the Bass--Quillen conjecture for vector bundles, it also reduces the proof of \Cref{intro-B-Q over val rings} to the case where $A$ is the coordinate ring of an affine open subset of some affine $R$-space.
\epp
For the sake of completeness, we present our Lindel-type result as follows.
We will only use its aforementioned special case when $Y=\{r=0\}$ for a nonunit $r\in A$ and $\#\mathbf{x}=1$ such that $A/rA\isoto B/rB$.
\bpropt [\cite{GL23}*{Proposition~4.4}]\label{variant of Lindel's lem}
Let $R$ be a semilocal Pr\"ufer domain. 
We fix
\begin{itemize}
\item [--]an irreducible, affine $R$-smooth scheme $X$ of pure relative dimension $d>0$, 
\item [--]a finitely presented closed subscheme $Y\subset X$ that is of $R$-fibrewise codimension $>0$, and 
\item [--]a finite subset $\mathbf{x} \subset X$.
\end{itemize}
If each fibre of $\mathbf{x}$ over a maximal ideal $\fm\subset R$ has fewer than $\max(\#\,k_{\mathfrak{m}},d)$ points, then there are 
\begin{itemize}
\item an affine open $W\subset X$ containing $\mathbf{x}$ and an affine open $U\subset \mathbb{A}_{R}^d$;
\item an \'etale $R$-morphism $f\colon W\to U$ fitting into the following Cartesian square 
\begin{equation*}
\qqq\begin{tikzcd}
W\cap Y \arrow[r, hook]  \arrow[d, equal]
 & W \arrow[d, "f"] \\
W\cap Y \arrow[r, hook]
& U, 
\end{tikzcd}
\x{\qqq where $f|_{W\cap Y}\colon W\cap Y \to U$ is a closed immersion.}
\end{equation*}
\end{itemize}
\epropt

\bremt\label{rem on Lindel's lem}
The assumption on $\# \mathbf{x}$ holds, for instance, if either $\mathbf{x}$ is a singleton or $d>\# \mathbf{x}$. 
Note that a certain assumption on $\# \mathbf{x}$ is necessary: when $X$ is a smooth affine curve over a finite field $\GGL=\bF_q$ and $\mathbf{x}\subset X(\bF_q)$, the resulting map $f$ from \Cref{variant of Lindel's lem} should give an injection $\mathbf{x} \hookrightarrow \mathbb{A}_{\bF_q}^1$, which is impossible as soon as $\# \mathbf{x} > q$.
\eremt

\section{Proof of the Bass--Quillen conjecture for torsors over valuation rings}
\label{sect-Bass-Quillen}
In this section, we prove the following result, the first main result of this paper.
\bthmt \label{B-Q over val rings}
Let $A$ be a ring that is ind-smooth over a Pr\"{u}fer ring $R$, and let $G$ be a totally isotropic reductive $R$-group scheme. Then, via pullback, we have
\[
\text{$H^1_{\mathrm{Nis}}(A,G)\isoto H^1_{\mathrm{Nis}}(\bA^N_A,G)$\q or, equivalently, \q $H^1_{\mathrm{Zar}}(A,G)\isoto H^1_{\mathrm{Zar}}(\bA^N_A,G)$}.
\]
\ethmt

\bpp[Sketch of proof of \Cref{Bass-Quillen for torsors} when $A$ contains a field] \label{pf-of-equi-char-BQ-torsors} 
Before proceeding further, we briefly outline a proof of the equi-characteristic case of \Cref{Bass-Quillen for torsors}, significantly simplifying Stavrova's argument in \cite{Sta22} and also avoiding reliance on Nisnevich's purity conjecture utilized in \cite{Ces22b}. 

Specifically, by employing Quillen's patching \Cref{Quillen-patching} and inducting on $N $, the problem reduces to showing that, for a totally isotropic reductive group scheme $ G $ over a regular local ring $ A $ containing a field, every generically trivial $ G $-torsor over $ \mathbb{A}_A^1$ is trivial. In other words, we need to show that the following composition map between pointed sets has a trivial kernel:
\[
H_{\et}^1(\mathbb{A}_A^1, G) \xrightarrow{\mu} H_{\et}^1(\mathbb{A}_K^1, G)  \xrightarrow{\nu}
H_{\et}^1(K(t), G),
\]
where $K=\Frac (A)$ and $t$ is the standard coordinate on the affine line. But this follows from:
\benum
  \item \label{ker(mu)-trivial}$\mu$ has trivial kernel. This follows from the relative form of the Grothendieck--Serre conjecture \cite{Fed22a}*{Theorem~1}: For a regular local ring $R$ containing a field $k$, with fraction field $K$, a totally isotropic reductive $R$-group scheme $G$, and an affine $k$-scheme $W$, no non-trivial $G$-torsor over $W\times_kR$ trivializes over $W\times_kK$. (It remains to take $W=\mathbb{A}_k^1$.)
  \item  \label{ker(nu)-trivial}$\nu$ has trivial kernel. By \cite{Gil02}*{corollaire~3.10}, a generically trivial $G$-torsor over $\mathbb{A}_K^1$ is isomorphic to $\lambda_*(\sO(1))$ for some cocharacter $\lambda$ of $G_K$, so it is itself trivial as $\sO(1)$ is so. \QED
\eenum
\epp

\bremt \label{mixed-char-rel-G-S}
The above argument would generalize to the mixed characteristic case if one could prove the relative form of the Grothendieck--Serre conjecture utilized in step \ref{ker(mu)-trivial} in this setting: For a DVR $D$, an essentially smooth local $D$-algebra $R$ with fraction field $K$, a totally isotropic reductive $R$-group scheme $G$, and a flat affine $D$-scheme $W$, no non-trivial $G$-torsor over $W\times_DR$ trivializes over $W\times_DK$.
\eremt

\bpp[Proof of \Cref{B-Q over val rings}]
   Since any section of $\mathbb{A}_A^N\to \Spec A$ induces sections to the pullbacks in \Cref{B-Q over val rings}, these pullbacks are injective. Thus, it remains to show that they are surjective. By Quillen patching \Cref{Quillen-patching}, we may assume that $R$ is a valuation ring by replacing $R$ with its localizations. 
   By a limit argument \cite{Gir71}*{VII, 2.1.6}, it suffices to assume that $A$ is $R$-smooth. 
   A relative limit argument involving \Cref{approxm semi-local Prufer ring} reduces us further to the case of a finite-rank $R$.

\textbf{Step 1:} \emph{$A$ is a polynomial ring over $R$.}

It suffices to show that every generically trivial $G$-torsor $\mathcal{E}$ over $R[t_1,\ldots,t_N]$ is trivial; \emph{a fortiori, } it descends to a $G$-torsor over $R$. We will argue by double induction on the pair $(N,\mathrm{rank}(R))$. If $N=0$, then, by convention $R[t_1,\ldots,t_N]=R$, the assertion follows from \cite{Guo24}*{Theorem~1.3} that a generically trivial $G$-torsor over a valuation ring is trivial. 
Now assume $N\geq 1$ and set $A\pr \ce R(t_N)[t_1,\ldots,t_{N-1}]$.
\epp
\bclt
The $G_{A\pr}$-torsor $\mathcal{E}_{A\pr}$ descends to a $G_{R(t_N)}$-torsor $\mathcal{E}_0$.
\eclt
\bpf[Proof of the claim]
Consider the natural projection  $\pi\colon \Spec R(t_N) \to \Spec R$. By definition, $R(t_N)$ is the localization of $R[t_N]$ with respect to the multiplicative system of monic polynomials. 
Thus, the closed fibre of $\pi$ is a singleton $\mathfrak{p}_{0}$. Furthermore, the local ring $R(t_N)_{\mathfrak{p}_0}$ is a valuation ring of $\mathrm{Frac}(R)(t_N)$, and its valuation restricts to the Gauss valuation on $R[t_N]$ associated to $R$  given by
\[ 
\tst R[t_N]\to \Gamma_R , \q  \sum_{i\ge 0} a_i t_N^i \mapsto \min_i v(a_i),
\]
where $v\colon R\to \Gamma_R$ is the (additive) valuation on $R$. In particular, we see that $R$ and $R(t_N)_{\mathfrak{p}_0}$ have the same value group. To apply the Quillen patching \Cref{Quillen-patching} and conclude, it suffices to show that the base change of $\cE_{A\pr}$ to $R(t_N)_{\mathfrak{p}}[t_1,\ldots,t_{N-1}]$ is trivial for every prime ideal $\fp \subset R(t_N)$.
 If $\mathfrak{p}=\mathfrak{p}_0$, then the above discussion implies that $\mathrm{rank}(R_{\pi(\fp)})=\mathrm{rank}(R)$, so
the desired triviality of the base change follows from induction hypothesis. If $\mathfrak{p}\neq \mathfrak{p}_0$, then $\pi(\mathfrak{p}) \in \Spec R$ is not the closed point, and we then have $\mathrm{rank}(R_{\pi(\fp)})<\mathrm{rank}(R)$. Therefore, by the induction hypothesis, $\cE_{R_{\pi(\mathfrak{p})}[t_1,\ldots,t_N]}$ is trivial, and hence its further base change along $R_{\pi(\mathfrak{p})}[t_1,\ldots,t_N]\to R(t_N)_{\mathfrak{p}}[t_1,\ldots,t_{N-1}]$ is also trivial.
\epf
  Since $\cE$ is generically trivial, so is $\cE_{A\pr}$. 
  Recall that the local ring $A\pr$ at the generic point of the closed fibre over $R$ is a valuation ring, by the Grothendieck--Serre \cite{Guo24}*{Theorem~1.3}, the torsor $\cE_{A\pr}$ is generically trivial on that closed fibre.
  Hence, by considering the pullback of $\mathcal{E}_{A\pr}$ along a general section $s\in \mathbb{A}_{R(t_N)}^{N-1}(R(t_N))$, we see that $\mathcal{E}_0$ is also generically trivial because it is trivial at the unique point lying over $\fm_R$. 
  By Lemma \ref{triviality over R(t)}, the torsor $\mathcal{E}_0$, and hence also $\mathcal{E}_{A\pr}$, is trivial. 
  Consequently, $\mathcal{E}$ is trivial away from the $R[t_1,\ldots, t_{N-1}]$-finite closed subset $\{f=0\} \subset \mathbb{A}_R^N=\Spec R[t_1,\ldots,t_N]$ for some monic polynomial $f\in R [t_N]$. 
  By \Cref{triviality over relative affine space for totally isotropic}, the $G$-torsor $\mathcal{E}$ must be trivial, which completes the induction process.

\textbf{Step 2:} \emph{$A$ is the localization $\wt{R}_S$ for a polynomial ring $\wt{R}\ce R[u_1,\ldots,u_d]$ with respect to some multiplicative subset $S\subset \wt{R}$.}

We wish to apply the `inverse' to Quillen patching \Cref{inverse patching}, with $\wt{R}$ as the base ring and $\wt{R}[t_1,\ldots,t_N]$ as the polynomial ring $A$. 
We still need to verify the assumptions there. First, by Step 1, any generically trivial $G$-torsor over $\wt{R}[t_1,\ldots,t_N]$ descends to a $G$-torsor over $\wt{R}$. 
Second, for any  generically trivial $G$-torsor $\mathcal{E}$ over $\wt{R}_S[t_1,\ldots,t_N]$, the restriction of $\mathcal{E}$ to each local ring of the 0-section
\[
   \Spec \wt{R}_S\simeq  \left\{t_1=\ldots=t_N=0\right\}\subset
 \mathbb{A}_{\wt{R}_S}^N
\]
is trivial (so extends to the trivial $G$-torsor over $\wt{R}$). 
Indeed, by Bass--Quillen in the field case, the restriction of $\cE$ to $\Frac (\wt{R}_S)[t_1,\ldots,t_N]$ is trivial. Thus, the restriction $\cE|_{\wt{R}_S}$ is generically trivial, and hence is Zariski-locally trivial (\Cref{G-S for constant reductive gps}). 
This verifies all the assumptions of \Cref{inverse patching}.

\textbf{Step 3:}  \emph{$A$ is an arbitrary smooth $R$-algebra.}
Let $\mathcal{E}$ be a generically trivial $G$-torsor over $\mathbb{A}_A^N$. 
Our goal is to show that $\mathcal{E}$ descends to a $G$-torsor over $A$.
Using Quillen patching \Cref{Quillen-patching}, we may assume that $A$ is an essentially smooth local algebra over a valuation ring $R$.
By localizing $R$, we can further assume that the ring map $R \to A$ is local.

We will argue by double induction on the pair $(\dim R, \dim A-\dim R)$ to show that $\mathcal{E}$ is even trivial. If $\dim R=0$, then $R$ is a field, so we conclude from the classical field case settled in \S~\Cref{pf-of-equi-char-BQ-torsors}. 
If $\dim A=\dim R$, then $A$ is also a valuation ring (\Cref{geom}~\ref{geo-iii}), so we reduce to the case already settled in Step 1. 
Assume now that $\dim R>0$ and $\dim A-\dim R>0$.

By \Cref{enlarge valuation rings}, we can enlarge $R$ (without changing $\dim R$) such that $A$ becomes the local ring $\sO_{X,x}$, where $X$ is an irreducible affine $R$-smooth scheme of pure relative dimension $d>0$, and $x\in X$ is a \emph{closed} point in the \emph{closed} $R$-fibre of $X$.
Note that $d=\dim A-\dim R$. 
By our Lindel-type result \S~\ref{geom_pres}, shrinking $X$ if needed, there are an \'etale $R$-morphism $f\colon X\to \mathbb{A}_R^d$ and a nonunit $r_0\in A_0\ce \sO_{\mathbb{A}_R^d,f(x)}$ such that
\[
   \text{$f$ \,\, induces a bijection \q \q
    $A_0/r_0A_0 \xrightarrow{\sim} A/r_0A$.}
\]
On the other hand, by induction hypothesis, $\mathcal{E}_{\mathbb{A}_{A_{\mathfrak{p}}}^N}$ is trivial for any $\fp\in \Spec A\backslash \{\mathfrak{m}_A\}$ (thus descends to the trivial $G$-torsor over $A_{\fp}$):
\begin{itemize}
  \item  either $\fp$ lies over a non-maximal ideal $ \mathfrak{q}\subset R$, in which case $R_{\mathfrak{q}} \to A_{\mathfrak{p}}$ is a local homomorphism,
      \[
      \dim R_{\mathfrak{q}} <\dim R, \q \text{ and } \q  \dim A_{\mathfrak{p}}-\dim R_{\mathfrak{q}}\le d=\dim A-\dim R;
      \]
  \item  or $R \to A_{\mathfrak{p}}$ is a local homomorphism, in which case
$$
\dim A_{\mathfrak{p}}-\dim R<\dim A-\dim R.
$$
\end{itemize}
Using Quillen patching \Cref{Quillen-patching}, we conclude that $\mathcal{E}_{\mathbb{A}_{A[1/{r_0}]}^N}$ descends to a $G$-torsor $\cF$ over $A[\frac{1}{r_0}]$. Since $\cF$ extends to a generically trivial $G$-torsor over $A$ (for example the restriction of $\mathcal{E}$ along any section $s\in \mathbb{A}_{A}^N(A)$), it must be trivial by \Cref{G-S for constant reductive gps}.
Now, \cite{Ces22a}*{Lemma~7.1} applies to the Cartesian square
\begin{equation*}
\begin{tikzcd}
\mathbb{A}_{A/r_0A}^N \arrow[hookrightarrow]{r} \arrow[d, "\sim"]
 & \mathbb{A}_A^N \arrow{d} \\
\mathbb{A}_{A_0/r_0A_0}^N \arrow[hookrightarrow]{r}
& \mathbb{A}_{A_0}^N,
\end{tikzcd}
\end{equation*}
we may glue $\mathcal{E}$ with the trivial $G$-torsor over $\mathbb{A}_{A_0[1/{r_0}]}^N $ to obtain a $G$-torsor $\cE_0$ over $\mathbb{A}_{A_0}^N$ that trivializes over $\mathbb{A}_{A_0[1/{r_0}]}^N$. 
By Step 2, the $G$-torsor $\cE_0$ is trivial, so $\mathcal{E}$ is trivial as well.
\QED

Building on \Cref{B-Q over val rings} and its equi-characteristic counterpart settled in \S~\ref{pf-of-equi-char-BQ-torsors}, we can now derive \Cref{Raghunathan-introd}, a generalization of Raghunathan's result, \cite{Rag89}*{Theorem~A}, to the relative context.
\bpf[Proof of \Cref{Raghunathan-introd}]
In view of \Cref{B-Q over val rings} and its equi-characteristic counterpart in full generality sketched in \S~\Cref{pf-of-equi-char-BQ-torsors}, it suffices to show the following. 
\bclt
If $(G,A)$ satisfies \Cref{Bass-Quillen for torsors}, then a $G$-torsor over $\bA^N_A$ is trivial if and only if it is trivial over $\bA^N_{K^s}$ and over $0_A\in \bA^N_A(A)$, where $K=\Frac A$ and $K^s$ is a separable closure of $K$.
\eclt
 
 \bpf[Proof of the claim]
  Let $\cE$ be a $G$-torsor over $\bA^N_A$ such that $\cE|_{\bA^N_{K^s}}$ and $\cE|_{0_A}$ are both trivial. First, note that it suffices to argue that $\cE$ is generically trivial so, in particular, we can assume that $A$ is a field (replacing $A$ by its fraction field): indeed, if so, then $\cE$ will descend to $A$ by \Cref{Bass-Quillen for torsors}, and will thus be trivial by checking along the zero section.
  \epf
 
 Assume now that $A=K$ is a field (Raghunathan's situation in \cite{Rag89}*{Theorem~A}). We assume by induction that $\cE|_{\bA^{N-1}_K}$ is already trivial,  where $\bA^{N-1}_K$ is regarded as a closed subscheme of $\bA^{N}_K$ via 
 \[
 \bA^{N-1}_K=\bA^{N-1}_K\times_K \{0\} \hookrightarrow \bA^{N-1}_K\times_K \bA^1_K=\bA^{N}_K.
 \]
Let $\pi\colon \bA^{N}_K \to \bA^{N-1}_K$ denote the projection onto the first $(N-1)$-coordinates, and write $\Spec L$ for the generic point of the target $\bA^{N-1}_K$. Then, taking generic fibre of $\pi$ we obtain a $G$-torsor $\cE|_{\bA^1_{L}}$ which is trivial over $\bA^1_{L^s}$ and over $0\in \bA^1_{L}(L)$. By a result of Raghunathan and Ramanathan (\emph{cf.,}~\cite{RR84}), this implies that $\cE|_{\bA^1_{L}}$ is trivial, so, in particular, $\cE$ is generically trivial, as desired.
 \epf

\section{$\bA^1$-homotopy theoretic affine representability in mixed characteristics}

\bpp[The $\infty$-category of motivic spaces]
\label{overview-Motivic-spaces}
Let $S$ be a quasi-compact quasi-separated scheme and let $\mathrm{Sm}_S$ denote the category of \emph{quasi-compact quasi-separated} (equivalently, finitely presented) $S$-smooth schemes\footnote{The quasi-compact quasi-separated assumption has two advantages: 1) the inclusion $\mathrm{Sm}_S \subset \Sch_{/S}$ is stable under taking pullbacks along smooth morphisms; 2) $\mathrm{Sm}_S$ is a small category, which cleanly avoids set-theoretical issues, although ultimately this is unnecessary (see \Cref{pass-between-Sm-Smaff}).}.  
The presentable $\infty$-category $\Spc(S)$ of motivic spaces over $S$ is constructed from the ordinary category $\mathrm{Sm}_S$ through the following steps:
\begin{itemize} 
\item [a)] Formally adjoin all small colimits in the $\infty$-categorical sense to create $\cP(\mathrm{Sm}_S)\ce \cP(\mathrm{Sm}_S, \cS)$, the $\infty$-category of presheaves of spaces on $\mathrm{Sm}_S$. 
\item [b)] Formally invert the class of all Nisnevich covering sieves to create $\Shv_{\mathrm{Nis}}(\mathrm{Sm}_S)$, the $\infty$-category of Nisnevich sheaves of spaces on $\mathrm{Sm}_S$. 
\item [c)] Formally contract the affine line $\bA^1_S$ by inverting the maps $\mathrm{pr}\colon \bA^1_U \to U$ for all quasi-compact quasi-separated (equivalently, all affine) $S$-smooth schemes $U$, to create $\Spc(S)$. 
\end{itemize}
\epp

Concretely, $\Shv_{\mathrm{Nis}}(\mathrm{Sm}_S) \subset \cP(\mathrm{Sm}_S)$ is the full subcategory spanned by those presheaves of spaces satisfying Nisnevich (\v{C}ech) descent. According to Morel--Voevodsky \cite{MV99}*{Proposition~1.4}, a presheaf satisfies Nisnevich descent if and only if it satisfies Nisnevich excision, that is, it sends the empty scheme $\emptyset$ to the final object $*\in \cS$ and sends Nisnevich squares to pullback squares in $\cS$ (see also \cite{SAG}*{Appendix~B.5.0.3}). As a result, the inclusion $\Shv_{\mathrm{Nis}}(\mathrm{Sm}_S) \subset \cP(\mathrm{Sm}_S)$ is stable under both limits and filtered colimits (since filtered colimits commute with pullbacks in $\cS$). Furthermore, since the pullback squares form a small set of conditions, this inclusion has a left exact (see \cite{HTT}*{Lemma~6.2.2.9}), accessible left adjoint (see \cite{HTT}*{Proposition~5.5.4.15}) 
\[
   L_{\mathrm{Nis}}\colon \cP(\mathrm{Sm}_S) \to \Shv_{\mathrm{Nis}}(\mathrm{Sm}_S),
\]
called the \emph{Nisnevich sheafification} functor. 
 Similarly, by \emph{loc.~cit.}, $\Spc(S)$ can be identified with the full subcategory of $\Shv_{\mathrm{Nis}}(\mathrm{Sm}_S)$ spanned by $\bA^1$-\emph{local objects}, i.e., those sheaves $\sF$ such that $\mathrm{pr}^*\colon \sF(U) \xrightarrow{\sim} \sF(\bA_U^1)$ for all $U \in \mathrm{Sm}_S$; moreover, the inclusion $\Spc(S) \subset \Shv_{\mathrm{Nis}}(\mathrm{Sm}_S)$ admits an accessible left adjoint (which is in general not left exact). 
Overall, $\Spc(S)$ is a presentable $\infty$-category, and the inclusion $\Spc(S) \subset \cP(\mathrm{Sm}_S)$ admits an accessible left adjoint
\[
L_{\mathrm{mot}}\colon \cP(\mathrm{Sm}_S) \to \Spc(S),
\]
called the \emph{motivic localization} functor. Moreover, by construction, the composition 
$$
\mathrm{Sm}_S \to \cP(\mathrm{Sm}_S) \xrightarrow{L_{\mathrm{mot}}} \Spc(S)
$$ is a functor such that postcomposing with it induces an equivalence
\begin{equation*}
  \Fun^{L}(\Spc(S),\cC) \xrightarrow{\sim} \Fun_{\Nis, \bA^1}(\mathrm{Sm}_S,\cC)
\end{equation*}
for any $\infty$-category $\cC$ with all small colimits, where $(-)_{\Nis, \bA^1}$ denotes the full subcategory of functors satisfying Nisnevich codescent and $\bA^1$-invariance, and $(-)^L$ denotes the full subcategory of colimit-preserving functors. Finally, by \cite{MV99}*{\S~2, Lemma~3.20}, $L_{\mathrm{mot}}$ can be described as
\begin{equation}\label{formula-Lmot}
\tst L_{\mathrm{mot}}=\left( L_{\mathrm{Nis}}\circ \Sing \right)^{\mathbf{N}}\ce \colim_{n\to \infty} \underbrace{L_{\mathrm{Nis}}\circ \Sing\circ \cdots \circ L_{\mathrm{Nis}}\circ \Sing}_{\x{$n$ times}}  .
\end{equation}
Here, the singular construction endofunctor $\Sing\colon \cP(\mathrm{Sm}_S)\to \cP(\mathrm{Sm}_S)$ is described by the formula 
\begin{equation}\label{Singular-constru}
  \Sing\sX\ce \colim_{\bDelta^{\mathrm{op}}}\sX(\Delta^{\bullet}_{\mathbf{Z}}\times_{\bbZ}-),
\end{equation}
where $\Delta^{\bullet}_{\mathbf{Z}}$ denotes the standard cosimplicial scheme over $\Spec \bbZ$. It is known that $\Sing$ takes values in the full subcategory $\cP_{\bA^1}(\mathrm{Sm}_S) \subset \cP(\mathrm{Sm}_S)$ of $\bA^1$-invariant presheaves and in fact it computes the left adjoint to the inclusion $\cP_{\bA^1}(\mathrm{Sm}_S) \subset \cP(\mathrm{Sm}_S)$. In particular, since the category $\bDelta^{\mathrm{op}}$ is sifted \cite{HTT}*{Lemma~5.5.8.4}, it follows that $L_{\mathrm{mot}}$ commutes with finite products.

\bremt \label{pass-between-Sm-Smaff}
In the construction of $\Spc(S)$ above, the category $\mathrm{Sm}_S$ can be replaced by either the category of all smooth $S$-schemes or just the full subcategory $\mathrm{Sm}_S^{\mathrm{aff}}$ of (absolutely) affine ones, without altering the resulting $\infty$-category  $ \Shv_{\mathrm{Nis}}(\mathrm{Sm}_S)$, and thus $\Spc(S)$, up to equivalences.  
For instance, to demonstrate that that the restriction functor
\begin{equation} \label{restr-Smaff-vs-Sm}
              \Shv_{\mathrm{Nis}}(\mathrm{Sm}_S) \to \Shv_{\mathrm{Nis}}(\mathrm{Sm}_S^{\mathrm{aff}})
\end{equation} 
is an equivalence, with its inverse given by the right Kan extension functor, we can apply \cite{HTT}*{Lemma~6.5.3.9}. This follows from the observation that every scheme admits a bounded\footnote{That is, if $(X_i)_{i=-1}^{\infty}$ denotes a hypercover, then there is an integer $N$ such that if $n>N$, then the canonical map $X_n\to \ker(\prod_{0\leq i\leq n}X_{n-1}\rightrightarrows \prod_{0\leq i<j\leq n}X_{n-2})$ is an isomorphism.} Zariski-hypercover where each term is affine: every scheme (resp., separated scheme) has an affine open cover such that every finite intersection is a separated scheme (resp., an affine scheme).
\eremt

\bpp[Representability results]\label{Represn-results}
Let $\sF\in \cP(\mathrm{Sm}_S)$ be a presheaf of spaces on $\mathrm{Sm}_S$. 
From the construction of $L_{\mathrm{mot}}(\sF)$ in (\ref{formula-Lmot}), there is a canonical morphism \begin{equation}\label{Sing-to-Motivic}
\Sing (\sF) \to L_{\mathrm{mot}}(\sF).
\end{equation} 
Following \cite{AHW18}*{Definition~2.1.1}, we say that $\sF$ is $\bA^1$-\emph{naive} if the restriction of (\ref{Sing-to-Motivic}) to $\cP(\mathrm{Sm}_S^{\mathrm{aff}})$ is an isomorphism. Under the equivalence (\ref{restr-Smaff-vs-Sm}), $\sF$ is $\bA^1$-naive if and only if $\Sing \sF \in \Spc(S)$. This follows, for instance, from the explicit formula (\ref{formula-Lmot}) for $L_{\mathrm{mot}}$. In the $\bA^1$-naive case, for every $X\in \mathrm{Sm}_S^{\mathrm{aff}}$, 
\[
\x{the map\q $\pi_0 \left( \Sing \sF\right)(X) \to \pi_0\left(L_{\mathrm{mot}}(\sF)(X)\right)\simeq [X, \sF]_{\bA^1}$\q  is bijective.}
\]

\epp

The following result is very useful for commuting geometric realizations and pullbacks.

\blemt
\label{commutes-geom-reali-pullback}
Let $\cC$ be an $\infty$-category. Consider the following Cartesian diagram in $\Fun(\bDelta^{\mathrm{op}},\cP(\cC))$:
\begin{equation*}
\begin{tikzcd}
(\sX_3)_{\bullet} \arrow[r] \arrow[d]
 & (\sX_2)_{\bullet} \arrow[d]\\
(\sX_1)_{\bullet} \arrow[r]
& (\sX_0)_{\bullet}.
\end{tikzcd}
\end{equation*}
If the simplicial object $(\pi_0\sX_0)_{\bullet} \in \Fun(\bDelta^{\mathrm{op}},\cP(\cC,\Set))$ is constant in the simplicial direction, then the following induced diagram of geometric realizations, valued in $\cP(\cC)$, is also Cartesian:
\begin{equation*}
\begin{tikzcd}
{|(\sX_3)_{\bullet}|} \arrow[r] \arrow[d]
 & {|(\sX_2)_{\bullet}|} \arrow[d]\\
{|(\sX_1)_{\bullet}|} \arrow[r]
& {|(\sX_0)_{\bullet}|}.
\end{tikzcd}
\end{equation*}
\elemt
\bpf  
Since the problem is sectionwise, we can formally replace $\cC$ with the category $\ast$, allowing us to work with $\cP(\ast) = \cS$. 
As argued in \cite{AHW17}*{Proof of Lemma~4.2.1}, the result then follows from the Bousfield--Friedlander theorem. Alternatively, we could cite \cite{Rez17}*{Proposition~5.4} or \cite{HA}*{Lemma~5.5.6.17}.  
\epf

\blemt[\emph{cf.}~\cite{AHW17}*{Theorem~5.1.3}]
\label{A1-naive-criterion}
Let $$\sF \in \Shv_{\mathrm{Nis}}(\mathrm{Sm}_S) \overset{(\ref{restr-Smaff-vs-Sm})}{\simeq}  \Shv_{\mathrm{Nis}}(\mathrm{Sm}_S^{\mathrm{aff}})$$ be a Nisnevich sheaf. If $\pi_0(\sF)$ is $\bA^1$-invariant on smooth affine $S$-schemes, then $\sF$ is $\bA^1$-naive. In particular, for every $U \in \mathrm{Sm}_S^{\mathrm{aff}}$, the natural map
\[
\tst \pi_0 \bigl( \Sing \sF\bigr)(U) \to [U, \sF]_{\bA^1}  \q  \text{ is bijective.}
\]
\elemt
\bpf
A key result proven in \cite{AHW17}*{Proposition~2.3.2} (which refines \cite{MV99}*{Proposition~1.4}) is that a presheaf $\sG$ on $\mathrm{Sm}_S^{\mathrm{aff}}$ is a Nisnevich sheaf if and only if it satisfies \emph{affine} Nisnevich excision, i.e., $\sG(\emptyset)=*$ and $\sG$ sends Nisnevich squares in $\mathrm{Sm}_S^{\mathrm{aff}}$ to pullback squares in $\cS$. 
Given a Nisnevich square 
\begin{equation*}
\begin{tikzcd}
V \arrow[hookrightarrow]{r} \arrow[d]
 & Y \arrow[d]\\
U\arrow[hookrightarrow]{r}
& X
\end{tikzcd}
\end{equation*}
in $\mathrm{Sm}_S^{\mathrm{aff}}$. 
Then we have the following Cartesian square in $\Fun(\bDelta^{\mathrm{op}},\cS)$:
\begin{equation*}
\begin{tikzcd}
\sF(X\times_{S}\Delta^{\bullet}_S) \arrow[r] \arrow[d]
 & \sF(U\times_{S}\Delta^{\bullet}_S) \arrow[d]\\
\sF(Y\times_{S}\Delta^{\bullet}_S)\arrow[r]
& \sF(V\times_{S}\Delta^{\bullet}_S).
\end{tikzcd}
\end{equation*}
By assumption, the simplicial set $\pi_0\sF(V\times_{S}\Delta^{\bullet}_S) $ is constant, we deduce from \Cref{commutes-geom-reali-pullback} (applied to $\cC=\ast $) that the induced diagram of geometric realizations is a Cartesian square in $\cS$. By the criterion mentioned above, $\Sing \sF$ is a Nisnevich sheaf on $\mathrm{Sm}_S^{\mathrm{aff}}$. The result then follows.
\epf

We can now finish the proof of Theorems \ref{affine-representability of H^1Nis} and \ref{naive-vs-A1-homotopy-G/H}.
\bpf
[Proof of \Cref{affine-representability of H^1Nis}]
Let $S$ be the spectrum of a Pr\"ufer ring.
Let $\bbB_{\Nis} G$ denote $L_{\mathrm{Nis}}(\bbB G)$, the Nisnevich sheafification of $\bbB G$. Then for every $U\in \mathrm{Sm}_S$, $(\bbB_{\Nis} G)(U)$ is (the nerve of) the groupoid of Nisnevich-locally trivial $G$-torsors over $U$, so that $\pi_0(\bbB_{\Nis} G)(U)=H^1_{\Nis}(U,G)$. By the Bass--Quillen \Cref{B-Q over val rings}, the presheaf $H^1_{\Nis}(-,G)$ on $\mathrm{Sm}_S^{\mathrm{aff}}$ is $\bA^1$-invariant, so \Cref{A1-naive-criterion} implies that $\Sing(\bbB_{\Nis}G)\simeq L_{\mathrm{mot}}(\bbB_{\Nis}G)$.
To conclude, note that since  $L_{\mathrm{mot}}(\bbB G)\simeq L_{\mathrm{mot}}(\bbB_{\Nis} G)$, by \Cref{A1-naive-criterion} we have:
\begin{align*}
  [U,\bbB G]_{\bA^1} & =\pi_0 \Map_{\Spc(S)}(L_{\mathrm{mot}}(U), L_{\mathrm{mot}}(\bbB G)) \\
   &  \simeq \pi_0\Map_{\Spc(S)}(L_{\mathrm{mot}}(U), L_{\mathrm{mot}}(\bbB_{\Nis} G))\\
   & =[U,\bbB_{\Nis} G]_{\bA^1}. \qedhere
\end{align*}
\epf

\bpf
[Proof of \Cref{naive-vs-A1-homotopy-G/H}]
Let $\ast =S\to \bbB_{\Nis} G$ denote the base point classifying the trivial $G$-torsor over $S$. We then have the following Cartesian diagram in $\Shv_{\mathrm{Nis}}(\mathrm{Sm}_S)$:
\begin{equation}\label{G/H-H-G-fibre-sequence}
\begin{tikzcd}
G/_{\Nis}H \arrow[r] \arrow[d]
 & \bbB_{\Nis} H \arrow[d]\\
\ast \arrow[r]
& \bbB_{\Nis} G,
\end{tikzcd}
\end{equation}
where the right vertical arrow is induced by the inclusion $H\subset G$ via functoriality, and $G/_{\Nis}H$ denotes the Nisnevich sheafification of $\mathrm{Sm}_S \ni U\mapsto G(U)/H(U)$. By assumption, we can identify $G/_{\Nis}H$ with the usual fpqc quotient $G/H$ restricted to $\mathrm{Sm}_S$. Evaluating the diagram (\ref{G/H-H-G-fibre-sequence}) at $\Delta_{S}^{\bullet}\times_S-$, we obtain a Cartesian diagram in $\Fun(\bDelta^{\mathrm{op}},\cP(\mathrm{Sm}_S^{\mathrm{aff}}))$. Since the reductive $S$-group $G$ is totally isotropic, \Cref{B-Q over val rings} implies that the simplicial object $$\pi_0(\bbB_{\Nis} G)(\Delta_{S}^{\bullet}\times_S-)=H^1_{\Nis}(\Delta_{S}^{\bullet}\times_S-,G)\in \Fun(\bDelta^{\mathrm{op}},\cP(\mathrm{Sm}_S^{\mathrm{aff}},\Set))$$ is constant in the simplicial direction. By taking geometric realizations, \Cref{commutes-geom-reali-pullback} yields the following Cartesian diagram in $\cP(\mathrm{Sm}_S^{\mathrm{aff}})$:
\begin{equation}\label{G/H-H-G-fibre-sequence}
\begin{tikzcd}
\Sing(G/H) \arrow[r] \arrow[d]
 & \Sing(\bbB_{\Nis} H) \arrow[d]\\
\ast \arrow[r]
& \Sing(\bbB_{\Nis} G).
\end{tikzcd}
\end{equation}
By the Bass--Quillen \Cref{B-Q over val rings} and the $\bA^1$-naive criterion \Cref{A1-naive-criterion}, we have  $\Sing(\bbB_{\Nis} H)$, $\Sing(\bbB_{\Nis} G) \in \Spc(S)$.
As the inclusion $\Spc(S) \subset \cP(\mathrm{Sm}_S^{\mathrm{aff}})$ is stable under all limits, it follows that $\Sing(G/H) \in \Spc(S)$, i.e., $G/H$ is $\bA^1$-naive. This completes the proof.
\epf


\end{document}